%%%%%%%%%%%% amsmacros.tex %%%%%%%%%%%%%%%%%%%%%%%%%%%%%%%%%%%%%%%%%%%%%%

\input amstex \input amssym.def \input amssym.tex \input amsppt.sty \nologo

%%%%%%%%%%%%%%%%%%%%%%%%%%%%%%%%%%%%%%%%%%%%%%%%%%%%%%%%%%%%%%%%%
%%%%%%%%%%%%%%%%%%%%%%%%%%%  Birkhauser macros %%%%%%%%%%%%%%%%%
%%%%%%%%%%%%%%%%%%%%%%%%%%%%%%%%%%%%%%%%%%%%%%%%%%%%%%%%%%%%%%%%%

\magnification=\magstephalf

\catcode`\@=12

\voffset=2.5truepc
\hoffset=2.5truepc
\hsize=5.15truein
\vsize=8.50truein %(including running head)
\def\halfhsize{2.5truein}

\tolerance=6000
\parskip=0pt
\parindent=18pt

\baselineskip=13.2pt

 \abovedisplayskip=12pt plus3pt minus2pt
 \belowdisplayskip=12pt plus3pt minus2pt
 \abovedisplayshortskip=12pt plus3pt minus2pt
 \belowdisplayshortskip=12pt plus3pt minus2pt

 \def\lskipamount{12pt}
 \def\lskip{\vskip\lskipamount plus3pt minus2pt}
 \def\lbreak{\par \ifdim\lastskip<\lskipamount
  \removelastskip \penalty-200 \lskip \fi}

 \def\lnobreak{\par \ifdim\lastskip<\lskipamount
  \removelastskip \penalty200 \lskip \fi}

% this is 14pt type for the title
\font\sectionfont=cmbx10 at 10pt
\font\autit=cmti12 %this is  type for the author's name
\font\titrm=cmbx12 at 14pt

%noindent first line of text
\def\ded#1{\rightline{{\it #1}}\medskip}

\def\sec#1{\goodbreak\vskip 1.5truepc\centerline{\hbox {{\sectionfont #1}}}
\nobreak\vskip1truepc\noindent\nobreak}

\def\refs{\vskip 1.5pc{\centerline {\bf References}}\vskip 4pt \noindent}
%then use \item macro \item{[1]}
%after references

%%%%%%%%%%%%%%%%%%%%%%%%%%%%%%%%%%%%%%%%%%%%%%%%%%%%%%%%%%%%%%%%%
%%%%%%%%%%%%%%%%%%% headers %%%%%%%%%%%%%%%%%%%%%%%%%%%%%%%%%
%%%%%%%%%%%%%%%%%%%%%%%%%%%%%%%%%%%%%%%%%%%%%%%%%%%%%%%%%%%%%%%%%

\newif\iftitlepage\titlepagetrue
\newtoks\runningtitle \runningtitle={}
\newtoks\runningauthor \runningauthor={}

\def\rhead{\iftitlepage{\hfil\global\titlepagefalse}\else
{\rm\hfil{\it \the\runningtitle}\hfill\folio}\fi}
\def\lhead{\rm\folio\hfill{\global\titlepagefalse \it\the\runningauthor}}

\headline={\ifodd\pageno\rhead\else\lhead\fi}

\footline={\hfil}
%\footline={\global\titlepagefalse\hfill}

%%%%%%%%%%%%%%%%%%%%%%%%%%%%%%%%%%%%%%%%%%%%%%%%%%%%%%%%%%%%%%%%%%%%
%%%%%%%%%%%%%%  additional macros (amsppt style) %%%%%%%%%%%%%%%
%%%%%%%%%%%%%%%%%%%%%%%%%%%%%%%%%%%%%%%%%%%%%%%%%%%%%%%%%%%%%%%%%%%%

\def\author{\bgroup\autit\bgroup}
\def\endauthor{\egroup\vskip 2truepc\egroup}

\def\title{\bgroup\bgroup\titrm\baselineskip20pt}
\def\endtitle{\egroup\vskip2truepc\egroup}

\def\Refs{\bgroup\refs\refstyle{A}\let\no=\key} \def\endRefs{\egroup}

\def\dedicatory{\bgroup\ded\bgroup}
\def\enddedicatory{\egroup\egroup}

\long\def\abstract#1\endabstract{
{\bf Abstract.} #1\par\vskip1.5pc}

\def\head#1\endhead{\sec{#1}}

\def\proclaim#1{\medskip{\noindent \bf  #1.}\bgroup\it}
\def\endproclaim{\egroup\medskip}

\def\demo#1{\medskip{\noindent\bf #1.}\rm\bgroup}
\def\enddemo{\egroup\medskip}

\def\remark#1{\medskip{\noindent\bf #1.}\rm\bgroup}
\def\endremark{\egroup\medskip}

%%%%%%%%%\font\ninerm=cmr9

\newbox\addressbox
\def\address{\setbox\addressbox=\vtop\bgroup
\hsize=\halfhsize\vglue3truepc\parindent=0pt
\obeylines}
\def\endaddress{\egroup} 

\newbox\addressboxtwo
\def\addresstwo{\setbox\addressboxtwo=\vtop\bgroup
\hsize=\halfhsize\vglue3truepc\parindent=0pt
\obeylines}
\def\endaddresstwo{\egroup}

%%%%%%%%%%%%%%%%%%%%%%%

%%%%%%%%%%%%%%%%%%%%%%%%%%

\title 
\centerline{Empirically Determined Ap\'ery -}
\centerline{Like Formulae for $\zeta(4n+3)$} 
\endtitle

\runningtitle{Formulae for $\zeta(4n+3)$}
\runningauthor{BORWEIN \& BRADLEY}

\author 
\centerline{Jonathan Borwein and David Bradley\footnote"\dag"{Research
supported by NSERC, the Natural Sciences and Engineering Research Council
of Canada.}}
\endauthor

%\ded{Dedicated to somebody}

\address 
Jonathan Borwein \& 
David Bradley,
Centre for Experimental
and Constructive Mathematics,
Simon Fraser University,
Burnaby B.C.,
Canada V5A 1S6.
jborwein\@cecm.sfu.ca
dbradley\@cecm.sfu.ca
\endaddress

%\loadeufm     % special fonts
%\input curly
%\TagsOnRight

\redefine\le{\leqslant}
\redefine\ge{\geqslant}
\redefine\i{\infty}

\define\({\left(}
\define\){\right)}
\define\[{\left[}
\define\]{\right]}
\define\il{\int_0^\i}

\define\Z{\text{\bf Z}}
\define\C{\text{\bf C}}
\define\a{\alpha}
\define\s{\sigma}
\define\z{\zeta}
\redefine\l{\lambda}
%\define\g{\gamma}
\define\G{\Gamma}
\define\Li#1{\operatorname{Li_#1}}

\interdisplaylinepenalty=500

\abstract
Some rapidly convergent formulae for special 
values of the Riemann Zeta function are given.  
We obtain a generating function formula for $\z(4n+3)$
which generalizes Ap\'ery's series for $\z(3)$, and appears
to give the best possible series relations of this type,
at least for $n<12$.  The formula reduces to a finite but
apparently non-trivial combinatorial identity.  The identity
is equivalent to an interesting new integral evaluation
for the central binomial coefficient.  We outline
a new technique for transforming and summing certain 
infinite series.  We also derive a beautiful formula which
provides strange evaluations of a large new class of non-terminating
hypergeometric series. 
It should be emphasized that our main results are 
shown equivalent but are still only conjectures.  
\endabstract

\head 1. Introduction \endhead
The Riemann Zeta function is 
$$
   \z(s) = \sum_{k=1}^\i \dfrac{1}{k^s},\qquad \Re(s)>1.\tag{1.1}
$$
Ap\'ery's formula \footnote{This formula is now commonly associated with
Ap\'ery's name because it featured in his proof of the irrationality of $\z(3).$
The formula actually goes back at least as far as \cite{5}.}
for $\z(3)$ is
$$
   \z(3) = \dfrac{5}{2}\sum_{k=1}^\i \dfrac{(-1)^{k+1}}{k^3 {2k\choose k}}.\tag{1.2}
$$
Extensive computation has suggested that there is no analogous formula
for $\z(5)$ or $\z(7)$.  In other words, if there exist relatively prime integers
$a$ and $b$ such that 
$$
   \z(5) = \dfrac{a}{b}\sum_{k=1}^\i \dfrac{(-1)^{k+1}}{k^5 {2k\choose k}},
$$
then $b$ is astronomically large.  Consider however, the following result of 
Koecher \cite{7}:
$$
   \z(5) = 2\sum_{k=1}^\i \dfrac{(-1)^{k+1}}{k^5{2k\choose k}} - \dfrac{5}{2}
            \sum_{k=1}^\i \dfrac{(-1)^{k+1}} {k^3{2k\choose k}} 
              \sum_{j=1}^{k-1}\dfrac{1}{j^2}.\tag{1.3}
$$
Inspired by this result, the second author, at the suggestion of the first author,
searched for additional
zeta identities of this sort using high-precision arithmetic and Maple's
lattice-based integer relations algorithms.  Within the hour, we were rewarded with
the following elegant new formula for $\z(7)$:
$$
   \z(7) = \dfrac{5}{2} \sum_{k=1}^\i \dfrac{(-1)^{k+1}}{k^7 {2k\choose k}}
          +\dfrac{25}{2} \sum_{k=1}^\i \dfrac{(-1)^{k+1}}{k^3 {2k\choose k}}
                   \sum_{j=1}^{k-1} \dfrac{1}{j^4}.\tag{1.4}
$$
Encouraged by this initial success, we searched for and found similar identities for $\z(9)$,
$\z(11)$, $\z(13)$, etc.  The representation for $\z(4n+3)$ has a convenient form in terms 
of a generating function (1.9), which is our main result (2.1).  It is curious that 
there is apparently no analogous generating function in the $4n+1$ case.  We refer the
reader to the discussion at the end of \S 8.  For now, it will be
advantageous to exhibit the recursive nature of the formulae in the $4n+3$ case.

We denote the power sum symmetric functions $P_r:=P_r^{(s)}(k)$ by $P_0:=1$ and
$$
   P_r^{(s)}(k) := \sum_{j=1}^{k-1} j^{-rs},\qquad r\ge 1.
$$
Next, we define a two-place function
$$
   \l(m,\prod_{j=1}^n P_{r_j}^{(s)}) 
:= \sum_{k=1}^\i\dfrac{(-1)^{k+1}}{k^m{2k\choose k}} \prod_{j=1}^n P_{r_j}^{(s)}(k).
$$
In (1.3) and \S 9, $s=2$ is relevant, but for now we are only interested in the case $s=4$. 
Therefore, to minimize symbol clutter we shall
occasionally repress the superscript, in which case $s=4$ should be assumed.
With Maple's help, the following list was produced:
$$
\align
  \dfrac{2}{5}\z(3) &= \l(3,P_0),\\
  \dfrac{2}{5}\z(7) &= \l(7,P_0) + 5\l(3,P_1),\\
  \dfrac{2}{5}\z(11) &= \l(11,P_0) + 5\l(7,P_1) -\dfrac{15}{2}\l(3,P_2) +
                         \dfrac{25}{2}\l(3,P_1^2),\\
  \dfrac{2}{5}\z(15) &= \l(15,P_0) + 5\l(11,P_1) -\dfrac{15}{2}\l(7,P_2) +
                         \dfrac{25}{2}\l(7,P_1^2)\\
                     &\qquad+ \dfrac{130}{6}\l(3,P_3) -\dfrac{225}{6}\l(3,P_1 P_2)
                      +\dfrac{125}{6}\l(3,P_1^3),
\tag{1.5}
\endalign
$$
etc.
The first list entry in (1.5) is just a restatement of Ap\'ery's formula (1.2) 
and the second list
entry is just a restatement of our formula (1.4).  From the list, it became clear to us that
the formula for $\z(4n+3)$ borrows the terms and coefficients from the formula 
for $\z(4n-1)$, except that the first argument of $\l$ is increased by $4$.  The
number of additional terms is equal to the number of partitions of $n$, and each
combination of power sum symmetric functions that occurs corresponds to a specific
partition of $n$.
Thus, we were led to conjecture that
$$
   \dfrac{2}{5}\z(4n+3) 
 = \sum_{j=0}^n \sum_{k=1}^\i\dfrac{(-1)^{k+1}}{k^{4j+3}{2k\choose k}}
   \sum_{\a\vdash n-j} c_{\a} P_{\a}^{(4)}(k),\tag{1.6}
$$
where the inner sum is over all partitions $\a$ of $n-j$, the $c_{\a}$ are rational numbers
to be determined, and if $\a=(\a_1,\a_2,\dots)$ is a partition of $n-j$, 
i.e. $\a_1+\a_2+\cdots = n-j$, then
$$
   P_{\a}^{(s)}(k) := \prod_{r\ge1}P_{\a_r}^{(s)}(k).
$$
Since it seemed plausible that a generating function could simplify matters, we rewrote our conjecture
(1.6) in the form
$$
\align
   \dfrac{2}{5}\sum_{n=0}^\i x^{4n}\z(4n+3)
 &= \sum_{n=0}^\i x^{4n} \sum_{j=0}^n \sum_{k=1}^\i \dfrac{(-1)^{k+1}}{k^{4j+3}{2k\choose k}}
     \sum_{\a\vdash n-j}c_{\a}P_{\a}^{(4)}(k)\\
 &= \sum_{j=0}^\i \sum_{n=j}^\i x^{4(n-j)}\sum_{k=1}^\i\dfrac{(-1)^{k+1}}{k^3{2k\choose k}}
      \(\dfrac{x}{k}\)^{4j}\sum_{\a\vdash n-j} c_{\a}P_{\a}^{(4)}(k)\\
 &= \sum_{s=0}^\i x^{4s}\sum_{k=1}^\i\dfrac{(-1)^{k+1}}{k^3{2k\choose k}} \dfrac{1}{1-x^4/k^4}
       \sum_{\a\vdash s}c_{\a}P_{\a}^{(4)}(k)\\
 &=  \sum_{k=1}^\i\dfrac{(-1)^{k+1}}{k^3{2k\choose k}}\dfrac{E_k(x^4)}{1-x^4/k^4},\tag{1.7}
\endalign
$$
where 
$$
   E_k(x) := \sum_{s=0}^\i x^{s} \sum_{\a\vdash s}c_{\a}P_{\a}^{(4)}(k).\tag{1.8}
$$
Let us denote the number of partitions of a non-negative integer $n$ by $p(n)$.
By convention, $p(0)=1$.
We suspected that $E_k(x)$ had a closed form which might be revealed by
determining enough of the coefficients in its power series.  Fortunately, due to
the recursive nature of the formulae we were able to extend the list (1.5) without
unduly straining Maple's lattice algorithms.  This was accomplished by introducing 
only $p(n)$ unknown coefficients for $\z(4n+3)$, rather than $\sum_{j=0}^n p(j)$, 
the actual number of terms involved.  Also, when the evidence warranted, we supplied
the coefficients of as many of the additional $p(n)$ terms as we confidently could,
based on our ability to recognize patterns and extrapolate from previously tabulated
values.
All identities so obtained were subsequently checked numerically -- typically to 250
significant digits.  

After having sufficiently extended the list (1.5), we were able to determine a good
many of the coefficients $c_{\a}$ for partitions $\a$ of small positive
integers, and hence the initial terms of the series expansion (1.8).
Maple's convert(series, ratpoly) feature then produced the following evaluations:
$$
\align
   E_1(x) &= 1,\\
   E_2(x) &= \dfrac{1+4x}{1-x},\\
   E_3(x) &= \dfrac{(1+4x)(16+4x)}{(1-x)(16-x)},\\
   E_4(x) &= \dfrac{(1+4x)(16+4x)(81+4x)}{(1-x)(16-x)(81-x)},\\
\endalign
$$
etc.  Thus we were led to conjecture that
$$
   E_k(x) = \prod_{j=1}^{k-1} \dfrac{j^4+4x}{j^4-x},
$$
and hence from (1.7) that
$$
   \sum_{n=0}^\i x^{4n}\z(4n+3)
 = \dfrac{5}{2}\sum_{k=1}^\i\dfrac{(-1)^{k+1}}{k^3{2k\choose k}}
   \dfrac{1}{1-x^4/k^4} \prod_{j=1}^{k-1} \dfrac{j^4+4x^4}{j^4-x^4},\qquad |x|<1.\tag{1.9}
$$
We restate (1.9) in the next section in the form of a conjectured
theorem, and discuss some of its
implications in the subsequent sections.

\head 2. A Generating Function Formula for $\z(4n+3)$\endhead
\proclaim{Theorem A (Conjectured)}
Let $z$ be a complex number.  Then
$$
   \sum_{k=1}^\i \dfrac{1}{k^3\(1-z^4/k^4\)}
 = \dfrac{5}{2} \sum_{k=1}^\i \dfrac{(-1)^{k+1}}{k^3{2k\choose k}}
                     \dfrac{1}{1-z^4/k^4}
                     \prod_{j=1}^{k-1} \dfrac{1+4 z^4/j^4}{1-z^4/j^4}.\tag{2.1}
$$
\endproclaim
\demo{Remark} 
Taking coefficients of $z^4$ in (2.1) yields our formula (1.4) for $\z(7)$.
Setting $z=0$ in (2.1) yields Ap\'ery's formula (1.2) for $\z(3)$.  
In general, taking
coefficients of $z^{4n}$ in (2.1) yields a rapidly convergent expansion for $\z(4n+3)$, the
$k$th term of which is a rational function of $k$ whose denominator is a power of $k$ 
times the central binomial coefficient, and whose
numerator is a symmetric function of partial harmonic sums in $1/j^4$.

More precisely, we denote the elementary symmetric functions by 
$$
   e_r^{(s)}(k) 
:= [t^r] \prod_{j=1}^{k-1}\(1+j^{-s}t\) 
 = \sum_{1\le j_1<j_2<\dots<j_r\le k-1} \(j_1 j_2\cdots j_r\)^{-s}
$$
and the complete monomial symmetric functions by
$$
   h_r^{(s)}(k)
:= [t^r] \prod_{j=1}^{k-1} \(1-j^{-s}t\)^{-1},
$$
where, as customary, $[t^r]$ means take the coefficient of $t^r$.
Then, assuming (2.1), we have the following
\proclaim{Corollary 1 (Equivalent to Conjectured Theorem A)} 
Let $n$ be a positive integer.  Then
$$
   \z(4n+3) 
 = \dfrac{5}{2} \sum_{j=0}^n \sum_{k=1}^\i \dfrac{(-1)^{k+1}}{k^{4j+3}{2k\choose k}}
                  \sum_{r=0}^{n-j} 4^r h_{n-j-r}^{(4)}(k) e_r^{(4)}(k).\tag{2.2}
$$
\endproclaim
\demo{Proof}
Extract the coefficient of $z^{4n}$ from each side of (2.1).  
% Apply $[z^n]$ to both sides of (2.1).
\enddemo

In light of the fact that both the complete symmetric functions and the elementary
symmetric functions can be expressed as rational linear combinations of the power sum
symmetric functions, it is possible to rewrite (2.2) in terms of the $P_{\a}$ of \S 1, as
in (1.6).
However, the formula for the coefficients $c_{\a}$
appears to be very complicated.  Thus, we have replaced the
sum over partitions in (1.6) with a much more manageable sum, at the expense
of introducing additional symmetric functions into the summand.

An additional consequence of (2.1) is an attractive formula which provides
strange evaluations for a large new
class of non-terminating hypergeometric series.
\proclaim{Corollary 2 (Equivalent to Conjectured Theorem A)} 
For all positive integers $n$, we have the formula
$$
\multline
   {}_6F_5 \bigg( \matrix n+1,\; n+1,\; 2n+in,\; 2n-in,\; in,\; -in\\
     n+1/2, n, 2n+1, n+1+in, n+1-in \endmatrix \bigg| -\tfrac{1}{4}\bigg) \\
 = \dfrac{2}{5}{2n\choose n}\prod_{j=1}^{n-1}\dfrac{n^4-j^4}{4n^4+j^4}.
\endmultline\tag{2.3}
$$
\endproclaim
\demo{Aside}  Throughout, we adhere to the standard notation
$$
   {}_pF_q \bigg( \matrix a_1, a_2,\dots, a_p \\
     b_1, b_2,\dots, b_q \endmatrix \bigg| z\bigg) := \sum_{k=0}^\i \dfrac{(a_1)_k(a_2)_k\cdots(a_p)_k}
         {(b_1)_k (b_2)_k\cdots (b_q)_k}\dfrac{z^k}{k!},
$$
where, as customary,
$$
   (a)_k := \dfrac{\G(a+k)}{\G(a)}= a(a+1)\cdots(a+k-1).
$$
\enddemo
\demo{Proof}
We can rewrite (2.1) as a formula for a non-terminating ${}_6F_5$:
$$
\multline
  \(\dfrac{1}{1-z^4}\)\\ \times {}_6F_5 \bigg(\matrix 2, 2, 1+z+iz, 1+z-iz, 1-z+iz, 1-z-iz\\
                   3/2,\;\;2+z,\;\;2-z,\;\;2+iz,\;\;2-iz\endmatrix \bigg| -\tfrac{1}{4}\bigg)\\
 = \dfrac{4}{5}\sum_{k=1}^\i \dfrac{1}{k^3\(1-z^4/k^4\)}.
\endmultline\tag{2.4}
$$
We note that both sides of (2.4) are meromorphic functions with simple poles at $z=\pm n$ and $z=\pm in$,
where $n$ is a positive integer.  We shall see that Corollary 2 is a consequence of equating residues
of both sides of (2.4) at the simple pole $z=n$.  If we denote the requisite residue by $R_n$, then from
the right side of (2.4), it is clear that
$$
   R_n = -\dfrac{1}{5n^2}.\tag{2.5}
$$
The residue calculation for the left side of (2.4) is more difficult.  We have
$$
\align
   R_n &= -\dfrac{1}{4z^3}\sum_{k=n-1}^\i \dfrac{(2)_k\, (2)_k\, (1\pm z\pm iz)_k}{(3/2)_k\, (1)_k\, (-4)^k}
              \prod^{k+1}\Sb {j=1}\\ j\ne n \endSb\dfrac{1}{j^4-z^4}\bigg|_{z=n}\\
       &= -\dfrac{n^4-z^4}{4z^3}\\
           &\quad\times\sum_{k=0}^\i \dfrac{(2)_{k+n-1}\, (2)_{k+n-1}\, (1\pm z\pm iz)_{k+n-1}}
              {(3/2)_{k+n-1}\, (1)_{k+n-1}\, (-4)^{k+n-1} } \prod_{j=1}^{k+n}\dfrac{1}{j^4-z^4}\bigg|_{z=n}\\
       &= -\dfrac{(n^4-z^4) (2)_{n-1}\,(2)_{n-1}\,(1\pm z\pm iz)_{n-1}}{4z^3 (3/2)_{n-1}\, (1)_{n-1}\, (-4)^{n-1}}
           \prod_{j=1}^n\dfrac{1}{j^4-z^4}\\
             &\quad\times\sum_{k=0}^\i \dfrac{(n+1)_k\,(n+1)_k\,(n\pm z\pm iz)_k}
           {(n+1/2)_k\,(n)_k\,(-4)^k}\prod_{j=n+1}^{n+k}\dfrac{1}{j^4-z^4}\bigg|_{z=n}\\
       &= -\dfrac{(2)_{n-1}\,(2)_{n-1}\,(1\pm n\pm in)_{n-1}}{n^3 (3/2)_{n-1}\,(1)_{n-1}\,4^n}\prod_{j=1}^{n-1}
           \dfrac{1}{n^4-j^4}\\
           &\quad \times {}_6F_5 \bigg( \matrix n+1,\; n+1,\; 2n+in,\; 2n-in,\; in,\; -in\\
     n+1/2, n, 2n+1, n+1+in, n+1-in \endmatrix \bigg| -\tfrac{1}{4}\bigg).\tag{2.6}\\
\endalign
$$
Comparing (2.5) and (2.6), it follows that
$$
\align
&{}_6F_5 \bigg( \matrix n+1,\; n+1,\; 2n+in,\; 2n-in,\; in,\; -in\\
     n+1/2, n, 2n+1, n+1+in, n+1-in \endmatrix \bigg| -\tfrac{1}{4}\bigg)\\
&=\dfrac{n 4^n (3/2)_{n-1}\,(1)_{n-1}}{5\, (2)_{n-1}\,(2)_{n-1}}\prod_{j=1}^{n-1}\dfrac{n^4-j^4}{4n^4+j^4}\\
&=\dfrac{2}{5}\(\dfrac{4^n (1/2)_n}{n!}\)\dfrac{n (1)_{n-1}}{(2)_{n-1}}\prod_{j=1}^{n-1}\dfrac{n^4-j^4}{4n^4+j^4}\\
&=\dfrac{2}{5}{2n \choose n}\prod_{j=1}^{n-1}\dfrac{n^4-j^4}{4n^4+j^4},
\endalign
$$
as required.  Thus we have shown that Corollary 2 follows from Conjectured Theorem
A.  That Corollary 2 implies Conjectured Theorem A now follows from 
Mittag-Leffler's Theorem.

\enddemo 
When $n=1$, the ${}_6F_5$ in Corollary 2 reduces to a ${}_4F_3$, and we obtain
\proclaim{Corollary 3}
$$
   {}_4F_3 \bigg( \matrix 2, 2, -i, i\\
     3/2,\, 1,\, 3\endmatrix \bigg| -\dfrac{1}{4}\bigg) = \dfrac{4}{5}.
$$
\endproclaim
\demo{Aside} Corollary 3 is true, and we have a proof.  However, since Corollary
3 is only a minor consequence of our conjectures, we delay the proof to the end of
\S 6, where the proof is used to illustrate some remarks we have to make
on our methods.
%\demo{Proof} Put $n=1$ in Corollary 2.
%\enddemo

\head 3.  Reduction to a Finite Identity \endhead
As we have said, (2.1) was originally a conjecture based on heavy experimental data.  
However, in the end,
we managed to reduce the problem to that of proving a finite combinatorial identity which
is beautiful in and of itself, and which we have, thus far, been unable to prove.
%, and which we have now proved.  
It is
$$
   \dfrac{5}{2}\sum_{k=1}^n {2k\choose k} \dfrac{k^2}{4n^4+k^4} \prod_{j=1}^{k-1}
              \dfrac{n^4-j^4}{4n^4+j^4} = \dfrac{1}{n^2},\qquad 1\le n \in \Z.\tag{3.1}
$$
The marvelous connection between the identity (3.1) and the 
conjectured generating function formula (2.1) is presented in the
%proof of the theorem below.
reduction below.
%\demo{Proof of Theorem}
\demo{Reduction}
By partial fractions we have for each positive integer $k$,
$$
   \dfrac{1}{1-z^4/k^4} \prod_{j=1}^{k-1} \dfrac{1+4z^4/j^4}{1-z^4/j^4}
 = \sum_{j=1}^k \dfrac{c_j(k)}{1-z^4/j^4},\tag{3.2}
$$
where
$$
   c_n(k) = \prod_{j=1}^{k-1} \(1+4n^4/j^4\) \bigg/ 
             \prod^k \Sb {j=1}\\ j\ne n \endSb \(1-n^4/j^4\),
   \qquad 1\le n \le k.\tag{3.3}
$$
Substituting (3.2) into the right hand side of (2.1) and interchanging
order of summation shows that (2.1) is equivalent to 
$$
   \sum_{k=1}^\i \dfrac{1}{k^3\(1-z^4/k^4\)}
 = \dfrac{5}{2} \sum_{j=1}^\i \dfrac{1}{1-z^4/j^4}
         \sum_{k=j}^\i \dfrac{(-1)^{k+1}c_j(k)}{k^3{2k\choose k}},
   \qquad z\in \C.
$$
Clearly, it suffices to prove that for all positive integers $n$,
$$
%   \sum_{k=n}^\i \dfrac{(-1)^{k+1}c_n(k)}{k^3{2k\choose k}}
   \sum_{k=n}^\i t_n(k)
 = \dfrac{1}{n^3},\tag{3.4}
$$
where 
$$
   t_n(k) := \dfrac{5(-1)^{k+1}c_n(k)}{2k^3{2k\choose k}}
   \qquad 1\le n \le k \in \Z.\tag{3.5}
$$
Our method of attack is to transform the infinite
sum (3.4) into a purely finite combinatorial identity.  This is accomplished
via analytic continuation of the summand combined with a process which 
might aptly be referred to as ``Gosper reflection''.  

Let $n$ be a fixed positive integer.  We wish to extend the definition
(3.5) to include values of $k$ less than $n$.  One approach is to convert
the products implicit in (3.5) into gamma functions.
Evidently
$$
\align
   c_n(k) 
 &= \lim_{x\to n} \dfrac{1-x^4/n^4}{1-x^4}
                  \prod_{j=1}^{k-1} \dfrac{1+4x^4/j^4}{1-x^4/(j+1)^4}\\
 &= \lim_{x\to n} \dfrac{\G(k\pm x\pm ix)}{\G(1\pm x\pm ix)}
                  \dfrac{\G^4(k+1)}{\G^4(k)} \dfrac{\G^4(n)}{\G^4(n+1)}
                  \dfrac{\G(n+1\pm x)}{\G(k+1\pm x)}\\
 &\qquad\qquad \times \dfrac{\G(n+1\pm ix)}{\G(k+1\pm ix)}
               \dfrac{\G(1\pm x)}{\G(n\pm x)} 
               \dfrac{\G(1\pm ix)}{\G(n\pm ix)}\\
 &= \dfrac{\G(k\pm n\pm in)}{\G(1\pm n \pm in)} \dfrac{k^4}{n^4}
    \dfrac{\G(2n+1)}{\G(n+k+1)} 
    \dfrac{1}{\G(k+1-n)} \dfrac{\G(n+1\pm in)}{\G(k+1\pm in)}\\
 &\qquad\qquad \times \dfrac{\G(n+1)}{\G(2n)} 
               \dfrac{(-1)^{n-1}}{(n-1)!} 
               \dfrac{\G(1\pm in)}{\G(n\pm in)}.\tag{3.6}\\
\endalign
$$
Here, and in the sequel, $\G(a\pm b\pm c)$ is shorthand for $\G(a+b+c)\G(a+b-c)\G(a-b+c)
\G(a-b-c)$.  It follows that for real $k$, one can define
$$
\align
   t_n(k) 
 &= -\dfrac{5}{2} \dfrac{e^{\pi i k}\G^2(k+1)}{\G(2k+1)}
    \dfrac{\G(k\pm n\pm in)}{\G(1\pm n\pm in)} \dfrac{k}{n^4} 
    \dfrac{\G(2n+1)}{\G(n+k+1)}\\
 &\times 
               \dfrac{1}{\G(k+1-n)}
               \dfrac{\G(n+1\pm in)}{\G(k+1\pm in)}
               \dfrac{\G(n+1)}{\G(2n)} \dfrac{(-1)^{n-1}}{(n-1)!}
               \dfrac{\G(1\pm in)}{\G(n\pm in)}.\tag{3.7}
\endalign
$$
Since $1/\G(k+1-n)=0$ when $k$ is an integer less than $n$, in view
of (3.4) it is necessary and sufficient to show that for all positive
integers $n$,
$$
   \sum_{k=0}^\i t_n(k)
 = \dfrac{1}{n^3}.\tag{3.8}
$$
To carry out the reflection process, we need to evaluate $t_n(k)$
when $k$ is a negative integer.  We shall see that when $k$ is a negative
integer, the rather forbidding-looking expression for $t_n(k)$
given by (3.7) takes a most attractive form.  From (3.7),
$$
\align 
   t_n(-1) 
 &= - \dfrac{5}{2n^4}\dfrac{(-1)^{n-1}}{(n-1)!}
      \dfrac{2n(n\pm in)(\pm in) n}{(\pm n\pm in)(-1\pm n\pm in)}
      \lim_{k\to-1}\dfrac{\G^2(k+1)}{\G(2k+1)\G(k+1-n)}\\
 &= \dfrac{5}{2n^2(1+4n^4)}\dfrac{(-1)^{n-1}}{(n-1)!}\\
 &\qquad \times
         \lim_{k\to-1}\dfrac{\G^2(k+2)}{(k+1)^2}
                 \dfrac{(2k+2)(2k+1)}{\G(2k+3)}
                 \dfrac{(k+1)}{\G(k+2)}
                 \prod_{j=0}^{n-1} (k-j)\\
 &= - \dfrac{5}{n(1+4n^4)}.\\
\endalign
$$
Let $j$ be a positive integer.  One can of course
evaluate $t_n(-j)$ directly from (3.7) by taking the limit as $k\to-j$,
just as we evaluated $t_n(-1)$ above.  However, it is preferrable to
introduce the following labour saving device.  For positive integer $k\ge n$,
define
$$
   \a_n(k) 
:= \dfrac{t_n(k)}{t_n(k+1)} 
 = \dfrac{-2k(2k+1)\((k+1)^4-n^4\)}{(k+1)^2(k^4+4n^4)}.
$$
For other values of $k$, define $\a_n(k)$ by the above expression on the far right-hand side.
Then for positive integer $k$,
$$
\align
   t_n(-k) 
 &= \a_n(-k)\a_n(1-k)\cdots\a_n(-2)t_n(-1)\\
 &= - \dfrac{5}{n(1+4n^4)}
      \prod_{j=2}^k \dfrac{2j(1-2j)}{(j-1)^2}
                    \dfrac{(j-1)^4-n^4}{4n^4+j^4}\\
 &= - \dfrac{5}{n(1+4n^4)}
      \dfrac{(2k)!}{2(k-1)!^2}\dfrac{n^4-1}{4n^4+k^4}
      \prod_{j=2}^{k-1}\dfrac{n^4-j^4}{4n^4+j^4}\\
 &= - \dfrac{5}{2n}{2k\choose k}\dfrac{k^2}{4n^4+k^4}
      \prod_{j=1}^{k-1}\dfrac{n^4-j^4}{4n^4+j^4}.\tag{3.9}
\endalign
$$
Set
$$
\align
   %p_n(k) &:= k^2(-1)^{n+1}(-k)_n(k+1\pm in)(k+2\pm in)\cdots(k+n-1\pm in),\\
   p_n(k) &:= k^2 \prod_{j=1}^{n-1} (k-j)(k+j+in)(k+j-in)\\
   %q_n(k) &:= (k-n\pm in),\\
   q_n(k) &:= (k-n+in)(k-n-in)\\
   r_n(k) &:= -2(k+n)(2k-1).
\endalign
$$
Then for all integers $k$, 
$$
   \dfrac{t_n(k+1)}{t_n(k)} 
 = \dfrac{p_n(k+1)}{p_n(k)} \dfrac{q_n(k)}{r_n(k+1)}.
$$
Furthermore, $q_n$ and $r_n$ share no linear factors differing by
an integer.  
According to Gosper's algorithm\footnote{We were led to consider Gosper's algorithm
when the first author attempted to get Maple to evaluate the sum (3.4) symbolically.
Mistyping `infinity' revealed that Maple could evaluate the resulting indefinite sum
for specific instances of the parameter $n$.} \cite{4}, there exists a polynomial $s_n$ of
degree no greater than $3n-3$ that satisfies
$$
   p_n(k) = s_n(k+1)q_n(k)-r_n(k)s_n(k)\tag{3.10}
$$
for all integers $k$.  Define
$$
   T_n(k) := \dfrac{r_n(k)s_n(k)t_n(k)}{p_n(k)}.\tag{3.11}
$$
Using (3.10), it is not hard to show 
that $T_n(k+1)-T_n(k)=t_n(k)$ for all integers $k$.
Note that since $t_n(-n)$ is finite and $p_n(-n)\ne 0= r_n(-n)$, we
have $T_n(-n)=0$.  It follows that 
$$
   T_n(m) = \sum_{j=-n}^{m-1} t_n(j),\qquad m-1\ge-n.
$$
Also, it is clear from (3.11) and (3.5) that 
$$
   \lim_{k\to\i} T_n(k) = 0.
$$
Thus (3.8) is equivalent to
$$
   T_n(0) = \sum_{j=-n}^{-1} t_n(j) = -\dfrac{1}{n^3}.\tag{3.12}
$$
Ideally, one would like to prove (3.12) using (3.11).  Unfortunately,
we do not know enough about the polynomials $s_n$ 
to infer the value $s_n(0)$ in general.  For specific values of $n$, 
we can use (3.10) to
solve for the unknown polynomial $s_n$ and hence, at least in 
principle, prove (3.12) for any
specific value of $n$.  However, using this approach to 
prove (3.12) in general would require an explicit
formula for the constant coefficient of the
possibly degree $3n-3$ polynomial $s_n$.  Of course, such a formula
can be inferred by assuming (3.12), but to us, at least, proving 
the formula directly seems a formidable task.  
However, substituting (3.9) into (3.12), it is readily apparent
that we need only prove the beautiful (and apparently non-trivial
\footnote{Professor Zeilberger \cite{10} has kindly informed us 
that (3.13) and its equivalent hypergeometric formulation (6.1)
do not fall under the scope of identities provable via the WZ method.
All attempts to prove (3.13) using WZ software, in particular, Zeilberger's
marvelous package `EKHAD' and Peter Paule's corresponding Mathematica implementation 
`zb\_alg.m', have failed.})
%and in particular the
%marvelous package `EKHAD', have failed.})
combinatorial identity
$$
   \dfrac{5}{2}\sum_{k=1}^n {2k\choose k} \dfrac{k^2}{4n^4+k^4}
                            \prod_{j=1}^{k-1} \dfrac{n^4-j^4}{4n^4+j^4}
 = \dfrac{1}{n^2}, \qquad n\in \Z^{+}.\tag{3.13}
$$
We discuss the identity (3.13) and some related results in the next section.
In \S 6, we examine the process of Gosper reflection in greater detail, where
it is revealed that identity
(3.13) and our conjectured generating function formula (2.1)
are in fact equivalent.

\head 4. A Combinatorial Identity \endhead
\proclaim{Lemma 4.1 (Equivalent to Conjectured Theorem A)} 
For all positive integers $n$, 
$$
   \sum_{k=1}^n \dfrac{5}{2}n^2k^2{2k\choose k} \dfrac{1}{4n^4+k^4}
                \prod_{j=1}^{k-1} \dfrac{n^4-j^4}{4n^4+j^4}
 = 1.
$$
\endproclaim
\noindent
Although we have verified Lemma 4.1 for all positive integers $n<300$, we have
so far been unable to find a proof.  The  
following equivalent proposition suggests one possible approach.
\proclaim{Proposition 4.2 (Equivalent to Conjectured Theorem A)}
For each positive integer $n$,
there exists an even polynomial $f_n$ of degree $2n$ such that
$$
   f_n(x)\prod_{j=1}^n\dfrac{x^2-j^2}{4x^4+j^4} 
 = 1-\dfrac{5}{2}\sum_{k=1}^n {2k\choose k}\dfrac{x^2 k^2}{4x^4+k^4}
   \prod_{j=1}^{k-1}\dfrac{x^4-j^4}{4x^4+j^4}.
$$
\endproclaim
\noindent
Clearly, Lemma 4.1 can be obtained from Proposition 4.2 if one sets $x=\pm 1,\pm 2,\pm 3,\dots \pm n$.
To see how we arrived at Proposition 4.2, let
$$
   \s_k(x) := \dfrac{5}{2}x^2 k^2{2k\choose k} 
              \prod_{j=1}^{k-1} (x^2+j^2)\tag{4.1}
$$
and define a sequence of functions $g_k$ recursively.  Put $g_0(x)=1$
for all $x$ and for $k>0$ let
$$
   g_{k-1}(x) - g_k(x) 
 = \dfrac{\s_k(x)}{4x^4+k^4} 
   \prod_{j=1}^{k-1} \dfrac{x^2-j^2}{4x^4+j^4}.\tag{4.2}
$$
Telescoping (4.2) would prove Lemma 4.1 if we could show that $g_n(n)=0$.
Define
$$
   f_k(x) := g_k(x) \prod_{j=1}^k\dfrac{4x^4+j^4}{x^2-j^2}.\tag{4.3}
$$
Then
$$
   g_k(x) = f_k(x) \prod_{j=1}^k\dfrac{x^2-j^2}{4x^4+j^4}.\tag{4.4}
$$
Clearly, $g_n(n)=0$ if $f_n(n)$ is finite.  In fact, the evidence strongly
suggests that each $f_k$ is a polynomial.  From (4.2) and (4.4)
it follows that
$$
   (4x^4+k^4)f_{k-1}(x) - (x^2-k^2)f_k(x) = \s_k(x), \qquad k>0.\tag{4.5}
$$
In particular, (4.2) and (4.5) imply that for all $x$,
$f_0(x)=1$, $f_1(x)=4x^2-1$, $f_2(x)=16x^4+4$, etc.  Now Proposition 4.2 is
obtained by telescoping (4.2) and writing $g_n$ in terms of $f_n$.

We remark that standard telescoping techniques prove the superficially similar
identity given in 
\proclaim{Proposition 4.3} For each positive integer $n$,
$$
   \dfrac{5}{4}\sum_{k=1}^n \dfrac{k^4 4^k}{4n^4+k^4}\prod_{j=1}^{k-1} \dfrac{n^4-j^4}
     {4n^4+j^4} = 1.
$$
\demo{Proof} 
Use 
$$
   \sum_{k=1}^n \(a_{k-1}-b_k\)\prod_{j=1}^{k-1} \dfrac{b_j}{a_j} = a_0 - b_n\prod_{j=1}^{n-1}
    \dfrac{b_j}{a_j}\tag{4.6}
$$
with 
$$
   a_k = \dfrac{4n^4+(k+1)^4}{4},\qquad b_k = n^4-k^4,\qquad k\ge 0.
$$
Standard telescoping proves (4.6) for any sequences of $a's$ and $b's$.  In our case, we
have
$$
   a_{k-1}-b_k = \dfrac{5}{4}k^4,\qquad b_n=0,
$$
and so
$$ 
   \dfrac{5}{4}\sum_{k=1}^n k^4 4^{k-1}\prod_{j=1}^{k-1}\dfrac{n^4-j^4}{4n^4+(j+1)^4}
 = \dfrac{4n^4+1}{4}.
$$
Now cross multiply and obtain
$$
   \dfrac{5}{4}\sum_{k=1}^n \dfrac{k^4 4^k}{4n^4+1}\prod_{j=1}^{k-1}\dfrac{n^4-j^4}{4n^4+(j+1)^4}
 = 1,
$$
from which the claimed identity easily follows.
\enddemo
If we try to play the same game using (4.6) to prove Lemma 4.1, it seems most natural to
define $ a_k := a_k(n) = 4n^4+(k+1)^4$ for $k\ge0$, and then choose $b_k:=b_k(n)$ so as to 
satisfy the recursion
$$
   (a_{k-1}-b_k)\prod_{j=1}^{k-1}\dfrac{b_j}{a_j} 
= \dfrac{5}{2}n^2 k^2{2k\choose k}\prod_{j=1}^{k-1}\dfrac{n^4-j^4}{4n^4+(j+1)^4},\qquad k\ge1.
\tag{4.7}
$$
If we can somehow show that $b_n(n)=0$, then (4.6) implies that
$$
   \sum_{k=1}^n \dfrac{5}{2} n^2 k^2 {2k\choose k} \prod_{j=1}^{k-1}\dfrac{n^4-j^4}{4n^4+(j+1)^4}
 = 4n^4+1,
$$
which, after cross multiplying, is easily seen to be equivalent to Lemma 4.1.
Now the recursion (4.7) is equivalent to 
$$
   \(a_{k-1}(n)-b_k(n)\)\prod_{j=1}^{k-1}\dfrac{b_j(n)}{n^4-j^4} 
 = \dfrac{5}{2}n^2 k^2{2k\choose k}.
$$
Thus, $b_n(n)=0$ is equivalent to
$$
   a_{n-1}(n)\prod_{j=1}^{n-1}\dfrac{b_j(n)}{n^4-j^4} = \dfrac{5}{2}n^4{2n\choose n}.
$$
i.e.
$$
   \prod_{j=1}^{n-1} \dfrac{b_j(n)}{n^4-j^4} = \dfrac{1}{2}{2n\choose n},
$$
which is an equivalent formulation of Lemma 4.1.

\head 5.  An Integral Identity \endhead
Here, we give an exquisite integral evaluation for the central binomial coefficient 
which is equivalent to Lemma 4.1 (3.13) and hence equivalent to our main conjecture.
\proclaim{Corollary 4 (Equivalent to Conjectured Theorem A)} 
For all positive integers $n$, we have the formula
$$
   \dfrac{1}{\pi}\il \dfrac{dy}{1+y^2}\prod_{j=0}^{n-1} \dfrac{4y^2-(j/n)^4}{y^2+(j/n)^4}
 = {2n\choose n}.
$$
\endproclaim

The proof that Corollary 4 is equivalent to Theorem A relies on
the following conjecture of Wenchang Chu (personal communication):  
\proclaim{Lemma 5.1 (Equivalent to Conjectured Theorem A)}
For all positive integers $n$, 
$$
   \sum_{k=1}^n \dfrac{2n^2}{k^2}\prod_{j=1}^{n-1}(j^4+4k^4)\bigg/{\prod^n \Sb
   j=1 \\ j\ne k \endSb (k^4-j^4)} 
 = {2n\choose n}.
$$
\endproclaim
\demo{Proof} We'll show that Lemma 5.1 and Lemma 4.1 are inverse pairs.
This fact is a special case of an inverse pair relationship given in \cite{2}
which is equivalent to 
$$
   f(n) = \sum_{k=r}^n \dfrac{a_n d_n + b_n c_n}{d_k} \dfrac{\phi(c_k/d_k;n)}
          {\psi_k(-c_k/d_k;n+1)}\,g(k)\tag{5.1}
$$
if and only if
$$
   g(n) = \sum_{k=r}^n \dfrac{\psi(-c_n/d_n;k)}{\phi(c_n/d_n;k+1)}\,f(k),\tag{5.2}
$$
where 
$$
   \phi(x;k) := \prod_{j=0}^{k-1} (a_j+x b_j),\quad 
   \psi(x;k) := \prod_{j=0}^{k-1} (c_j+x d_j),
$$
and
$$
   \psi_m(x;k) := \prod^{k-1} \Sb j=0 \\ j\ne m \endSb (c_j+x d_j).
$$
Setting $r=1$, $a_j=j^4$, $b_j=4$, $c_j=j^4$, $d_j=1$, $f(k) = 10 k^2 {2k\choose k}(-1)^k$,
$g(n)=1/n^2$ in the inverse pair (5.1), (5.2) yields the claimed inverse pair relationship
between Lemma 4.1 and Lemma 5.1.
\enddemo

We now proceed to show that Corollary 4 is equivalent to Lemma 5.1.
By a suitable change of variable, the integral identity in Corollary 4
can be rewritten in the form
$$
   \dfrac{4n^2}{\pi} \il dx \prod_{j=1}^{n-1} (4x^2-j^4)\bigg/ \prod_{j=1}^n (x^2+j^4)
 = {2n\choose n}.\tag{5.3}
$$
In view of the partial fraction expansion (3.2), we can rewrite the
integrand of (5.3), obtaining the equivalent identity
$$
\align
{2n\choose n}
& = \dfrac{4n^2}{\pi} \il (-1)^{n+1} \sum_{k=1}^n \dfrac{k^4}{n^4} \dfrac{c_k(n)}{k^4+x^2}\,dx\\
&= (-1)^{n+1}\sum_{k=1}^n \dfrac{2k^2}{n^2}c_k(n),
\endalign
$$
which, in view of the definition (3.3) of the numbers $c_k(n)$, is precisely the statement
of Lemma 5.1.

\head 6.  Some Remarks on Reflection \endhead
We can rewrite Lemma 4.1 or (3.13) in hypergeometric notation as
$$
\multline
   {}_6F_5 \bigg( \matrix 2,\;\; 3/2,\;\; 1+n,\;\; 1-n,\;\; 1+in,\;\; 1-in\\
     1,\, 2+n+in,\, 2+n-in,\, 2-n+in,\, 2-n-in \endmatrix \bigg| -4\bigg) \\
 = \dfrac{4n^4+1}{5n^2},
\endmultline\tag{6.1}
$$
an apparently new \footnote{See the remarks at the end of this section.}
strange evaluation of a terminating ${}_6F_5.$
We can also rewrite (2.1) as a formula for a non-terminating ${}_6F_5$:
$$
\multline
   {}_6F_5 \bigg(\matrix 2,\,2,\,1+z+iz,\,1+z-iz,\,1-z+iz,\,1-z-iz\\
                   3/2,\;\;2+z,\;\;2-z,\;\;2+iz,\;\;2-iz\endmatrix \bigg| -\tfrac{1}{4}\bigg)\\
 = \dfrac{4}{5}\sum_{k=1}^\i \dfrac{1-z^4}{k^3\(1-z^4/k^4\)}.
\endmultline\tag{6.2}
$$
Observe the dual nature of (6.1) and (6.2).  Our process of Gosper reflection has
taken a non-terminating ${}_6F_5$ at $-1/4$, and transformed it into a terminating
${}_6F_5$ at $-4$, in which certain of the numerator parameters and denominator parameters have
been exchanged and shifted.

We can see the dual results of reflection in another way.
Let $z^4=-n^4/4$ in (2.1).  The right-hand side terminates, yielding
$$
   \sum_{k=1}^\i \dfrac{4k}{4k^4+n^4} 
=  \dfrac{5}{2}\sum_{k=1}^n \dfrac{4^k}{{2k\choose k}}\dfrac{k}{n^4+4k^4}
     \prod_{j=1}^{k-1} \dfrac{n^4-j^4}{n^4+4j^4}.\tag{6.3}
$$
On the other hand, standard techniques show that 
$$
\align
   \sum_{k=1}^\i \dfrac{4k}{4k^4+n^4} 
&= \dfrac{1}{2in^2}\bigg\{\psi\(1-in\(\tfrac{1+i}{2}\)\)+\psi\(1+in\(\tfrac{1+i}{2}\)\)
      -\psi\(1+n\(\tfrac{1+i}{2}\)\)\\
&\qquad -\psi\(1-n\(\tfrac{1+i}{2}\)\)\bigg\}\\
&= \dfrac{1}{2n}\sum_{k=1}^n \dfrac{1}{\(k-n/2\)^2 + n^2/4}.\tag{6.4}
\endalign
$$
Comparing (6.3) and (6.4) yields the following identity:
$$
   \dfrac{5}{2}\sum_{k=1}^n \dfrac{4^k}{{2k\choose k}}\dfrac{k}{n^4+4k^4}
     \prod_{j=1}^{k-1} \dfrac{n^4-j^4}{n^4+4j^4}
 = \dfrac{1}{2n}\sum_{k=1}^n \dfrac{1}{\(k-n/2\)^2 + n^2/4}.\tag{6.5}
$$
Now compare the left-hand sides of (6.5) and (3.13).

The astute reader will observe a close relationship between the right side of (2.3) in Corollary 2
and the summand of Lemma 4.1.  In fact, Gosper reflection applied to Corollary 2 yields the identity
in Lemma 4.1.  Since the proof of this mirrors the development of \S 3, we omit the details.  We
remark however, that Gosper reflection easily proves any specific instance of Corollary 2.  For the
sake of brevity, we illustrate this assertion in the case $n=1$ i.e. Corollary 3.
Writing $t(k)$ for the summand of
Corollary 3, we have
$$
   t(k) = \dfrac{(2)_k (2)_k (\pm i)_k (-1/4)^k }{(3/2)_k (1)_k (1)_k (3)_k}
        = \dfrac{(k+1)^2 \G(k\pm i) \G(1/2) (-1/4)^k}{\G(\pm i)\G(3/2+k)\G(k+3)}.
$$
It follows that $t(-1)=0$ and $t(-k)=0$ for $3\le k\in \Z$.  Since
$$
   \sum_{k=-\i}^\i t(k) = 0,
$$
it follows that
$$ 
   \sum_{k=0}^\i t(k) = -t(-2) = \dfrac{4}{5},
$$
which proves Corollary 3.

To conclude this section, we'd like to offer evidence in support of
our claim that evaluations
(6.1) and Corollary 2 are indeed new.  
After surveying the standard references (eg. \cite{3,1,9}) in
the vast hypergeometric literature and consulting many
of the experts in this area, we have been unable to uncover
anything remotely like (6.1) or Corollary 2. 
Hypergeometric summations in
which the main argument is different from 1 are rare enough.
Exceedingly rare are summations with complex parameters such as in
(6.1) or Corollary 2, and neither of our evaluations appears
to have a natural generalization.  For example, there appears
to be no generalization of either formula in which $in$ is replaced by
a general parameter $m$ for example.

\head 7. Algorithms and Complexity \endhead
The formulae developed herein lend themselves quite easily to numerical computation.
Algorithms based on Ap\'ery's formula (1.2), Koecher's formula (1.3), and Bradley's
formula (1.4) are particularly simple and are given below.

$$
\align
&\text{Zeta3}:=\text{proc}(d)\\
&\# \text{Compute } \z(3) \text{ to } d \text{ digits using }(1.2).\\
&N:=1+\lfloor 5d/3\rfloor;\, \text{Digits}:=d;\, c:=2;\,s:=0;\\
&\text{for } n \text{ from }1 \text{ to } N \text{ do}\\
&\qquad    s:=s+(-1)^{n+1}/( n^3 c);\; c:=c(4n+2)/(n+1);  \\    
&\text{od};\\
&\text{Return}(5s/2);\\
\\
&\text{Zeta5}:=\text{proc}(d)\\
&\# \text{Compute } \z(5) \text{ to } d \text{ digits using }(1.3).\\
&N:=1+\lfloor 5d/3\rfloor;\, \text{Digits}:=d;\, a:=0;\,c:=2;\,s:=0;\\
&\text{for } n \text{ from } 1 \text{ to } N \text{ do }\\
&\qquad   g:=1/n^2;\; s:=s+(-1)^{n+1}(4g-5a)/(n^3 c);\\
&\qquad   c:=c(4n+2)/(n+1);\; a:=a+g;\\
&\text{od}; \\
&\text{Return}(s/2);\\
\\
&\text{Zeta7}:=\text{proc}(d)\\
&\# \text{Compute }\z(7) \text{ to } d \text{ digits using }(1.4).\\
&N:=1+\lfloor 5d/3\rfloor;\, \text{Digits}:=d;\, a:=0;\,b:=0;\,c:=2;\,s:=0;\\
&\text{for } n \text{ from }1 \text{ to } N \text{ do}\\
&\qquad	g:=1/n^4;\; s:=s+(-1)^{n+1}(5a+g)/(n^3 c);\\
&\qquad c:=c(4n+2)/(n+1);\; a:=a+g;\\
&\text{od};\\
&\text{Return}(5s/2); \\
\endalign
$$
By Stirling's asymptotic formula for the gamma function, it readily follows that
$$
   {2k\choose k} \sim \dfrac{4^k}{\sqrt{\pi k}},\qquad k\to\i,
$$
and thus all formulae we have discussed yield $2$ binary digits per term asymptotically.
Since $\log 4/\log10 \approx 3/5$, this amounts to slightly better than $1.2$ decimal
digits per term.
This should be contrasted with the definition (1.1) which is asymptotically useless,
yielding $0$ digits per term.  For example, computing $\z(3)$ from the definition (1.1)
and applying the integral test to the tail of the series shows that the $n$th
tail drops off like $O(1/n^2)$.  Thus each successive digit requires computing
$\sqrt{10}$ times as many terms as its predecessor.  To get $d$ digits, $O(10^{d/2})$
operations are involved.  On the other hand, it's not hard to see that the algorithms we
have presented require only $O(d)$ operations to compute $d$ digits.

In the preceding discussion, we spoke of operations as operations (such as multiplication
or division) on numbers.  A more realistic evaluation of run times must
take into account operations on digits.  If we take as given that the cost of multiplying two
$d$ digit numbers is $O(d(\log d)\log \log d)$, then a crude upper bound on the run time for
computing $d$ digits using our Ap\'ery-like algorithms is
$$
   O\big(\sum_{j=1}^d j(\log j) \log\log j\big) = O(d^2(\log d)\log \log d).
$$
However, it is possible to adapt these algorithms using the method of Karatsuba \cite{6}
to yield the highly respectable run time $O(d (\log d)^3 \log \log d).$
We coded our Ap\'ery-like algorithms (sans Karatsuba's optimization) in Maple V 
(release 3) and ran them on an Indy
R4600PC 100 MHz Silicon Graphics Workstation.
The following table compares the run times in cpu seconds with Maple's built-in implementation
of the Riemann zeta function.
\medskip \medskip

\hbox to\hsize{\hfil\bf Table of Run Times\hfil}\nobreak\bigskip
\def\mspace{&&&&\cr\noalign{\vskip-5pt}}
\def\hmspace{\hskip6pt}
\vglue -.3in
$$
\vbox{
\offinterlineskip
\halign{
\strut\vrule\hfil\hmspace$#$\hmspace\vrule&&\hfil\hmspace$#$\hmspace\vrule\cr
\noalign{\hrule}
\mspace
&\z(3)\,\,\,\, &\z(5)\,\,\,\, &\z(7)\,\,\,\,
&\text{Digits}\cr
\mspace
\noalign{\hrule}
\mspace
\text{Ap\'ery-like}  & 0.4561  & 1.8720  & 2.8141  & 200 \cr 
\mspace \noalign{\hrule}\mspace
\text{Maple} & 8.1720  & 8.4600  & 8.3462 & 200 \cr 
\mspace \noalign{\hrule}\mspace
\text{Ap\'ery-like}  & 1.1401  & 5.5019  & 8.0399  & 300 \cr 
\mspace \noalign{\hrule}\mspace
\text{Maple} & 28.0742  & 28.1819  & 28.3860   & 300\cr 
& & & & \cr
\noalign{\hrule}
}
}
$$

\head 8. Other Dirichlet Series \endhead
For all positive integers $n$ and all real $k$, let
$$
   d_n(k) := \dfrac{5n^3}{2k^3}c_n(k).
$$
Then (3.4) becomes
$$
   \sum_{k=n}^\i \dfrac{(-1)^{k+1}}{{2k\choose k}}d_n(k) = 1, \qquad 1\le n\in\Z.
$$
Thus, for any sequence $a_1,a_2,\dots$, we may write
$$
   \dfrac{a_n}{n^s} 
 = \dfrac{a_n}{n^s}\sum_{k=n}^\i \dfrac{(-1)^{k+1}}{{2k\choose k}}d_n(k).\tag{8.1}
$$
Let us suppose that $\sum n^{-s} a_n$ is absolutely convergent.
Summing (8.1) on $n$ and interchanging the order of summation, we get 
$$
   \sum_{n=1}^\i \dfrac{a_n}{n^s} 
 = \sum_{k=1}^\i \dfrac{(-1)^{k+1}}{{2k\choose k}}\sum_{j=1}^k \dfrac{a_j}{j^s} d_j(k).\tag{8.2}
$$
This gives a ``formula'' for any absolutely convergent Dirichlet series.  However, (8.2) does
not appear to be of much use, except in special cases where we can take advantage of known
properties of the numbers $d_n(k)$.  For example, since
$$
   \sum_{j=1}^k c_j(k) = 1
$$ 
for all $1\le j\le k\in\Z$, putting $s=3$ and $a_1=a_2=\dots=1$ in (8.2) recovers Ap\'ery's
formula (1.2).  

Unfortunately, there seems to be no way to make use of (8.2) or the ideas
of \S3 to obtain a generating function analogue of our result (2.1) for $\z(4n+1)$.  
Since (2.1) started with Ap\'ery's formula (1.2) for $\z(3)$, one might expect
that a generating function analogue of
(2.1) for $\z(4n+1)$ would be based on Koecher's formula (1.3) for $\z(5)$
and derive from recurrence properties akin to those implicit in the list (1.5).
However, none of the formulae for $\z(9)$ that we discovered (and we have good
reason to believe there are no others) bears the necessary relationship to (1.3).

We should also point out that even in the $4n+3$ case, much work remains to be done, as
there are several Ap\'ery-like formulae for $\z(7)$, $\z(11)$, etc. which do not arise from 
our generating function (2.1).  In the $4n+1$ case, the proliferation of formulae appears
to be even greater.  We have 
created code for systematically listing the formulae for $\z(13)$, and ran
the code for two months or so.  The resulting file is over three thousand lines long and
contains hundreds and hundreds of independent formulae, all having the characteristic 
power of $k$ and
central binomial coefficient in the denominator, accompanied by harmonic-like sums in
the numerator. 
Classifying the myriad relations and interrelations
amongst these sums for the various even/odd zeta values would be a huge project indeed.

\head 9.  Addendum \endhead
As we later learned, Koecher \cite{8} had given a very simple proof of the following 
generating function for $\z(2n+1)$, namely
$$
   \sum_{k=1}^\i \dfrac{1}{k^3\(1-z^2/k^2\)}
 = \sum_{k=1}^\i \dfrac{(-1)^{k+1}}{k^3{2k\choose k}}\bigg(\dfrac{1}{2}+\dfrac{2}{1-z^2/k^2}\bigg)
      \prod_{j=1}^{k-1}\(1-z^2/j^2\).\tag{9.1}
$$
If $n$ is a non-negative integer, 
%applying $[z^n]$ to both sides of (9.1) produces 
extracting the coefficient of $z^n$ from each side of (9.1) produces
the formula
$$
\multline
   \z(2n+3) = \dfrac{5}{2}\sum_{k=1}^\i \dfrac{(-1)^{k+1}}{k^3{2k\choose k}}(-1)^n e_n^{(2)}(k)\\
            + 2\sum_{j=1}^n \sum_{k=1}^\i \dfrac{(-1)^{k+1}}{k^{2j+3}{2k\choose k}}(-1)^{n-j}
                       e_{n-j}^{(2)}(k),
\endmultline\tag{9.2}
$$
from which (1.2) and (1.3) follow as special cases.  Here, the $e_r^{(s)}(k)$ are the 
elementary symmetric functions defined in \S 2. 

Despite the fact that Koecher's generating
function (9.1) gives formulae for all odd Zeta values, there is a very real sense in which 
(9.1) is inferior to our generating function (2.1).  
In (9.1), among other things, the fourth powers that feature in (2.1) are replaced 
by squares.  This results in redundant terms in his zeta formula (9.2) for $n>1$.  
For example, $n=2$ in (9.2) yields 
$$
\multline
   \z(7) = 2\sum_{k=1}^\i \dfrac{(-1)^{k+1}}{k^7{2k\choose k}} - 2\sum_{k=1}^\i \dfrac{(-1)^{k+1}}
                      {k^5{2k\choose k}} \sum_{j=1}^{k-1}\dfrac{1}{j^2}\\
           +\dfrac{5}{2}\sum_{k=1}^\i \dfrac{(-1)^{k+1}}{k^3{2k\choose k}}\sum_{1\le j<l\le k-1}
             \dfrac{1}{j^2 l^2},
\endmultline\tag{9.3}
$$
which should be compared with our more compact formula (1.4).  
To enable a more detailed comparison, we rewrite (9.3) in the notation of \S1.  Then (9.3)
becomes
$$
   \z(7) = 2\l(7,P_0^{(2)}) - 2\l(5,P_1^{(2)}) + \dfrac{5}{4}\l(3,P_1^{(2)}P_1^{(2)})
            -\dfrac{5}{4}\l(3,P_2^{(2)}),\tag{9.4} 
$$
whereas (1.4) is simply
$$
   \z(7) = \dfrac{5}{2}\l(7,P_0^{(4)}) +\dfrac{25}{2}\l(3,P_1^{(4)}).\tag{9.5}
$$
Since $P_0^{(4)}=P_0^{(2)}=1$ and $P_1^{(4)}=P_2^{(2)}$, we see that the middle two
terms of (9.4) are redundant.  Indeed, lattice-based reduction shows that 
$$
   2\l(7,P_0^{(2)}) + 8\l(5,P_1^{(2)}) -5\l(3,P_1^{(2)}P_1^{(2)}) +55\l(3,P_2^{(2)}) = 0.
$$
As far as we can tell, in contrast with the formulae derived from (9.2),
there are no redundant terms in our formulae for $\z(4n+3)$ 
which come from (2.1), at least for $n<12.$  It goes without saying that, despite our
best efforts, Koecher's proof of (9.1) apparently cannot be adapted to prove (2.1).
It seems that (1.4), and more 
generally (2.1), is a much deeper result.  We should also point out that
merely bisecting Koecher's generating function (9.1) will not yield (2.1), nor any new zeta
formulae.

\copy\addressbox
\bye